# Generalized-hypergeometric solutions of the biconfluent Heun equation

D.Yu. Melikdzhanian[1,2] and A.M. Ishkhanyan[1,3]

[1]Russian-Armenian University, Yerevan, 0051 Armenia
[2]Vekua Institute of Applied Mathematics, Tbilisi State University, Tbilisi, 0186 Georgia
E-mail: davidmelikdzhanian@gmail.com
[3]Institute for Physical Research, Ashtarak, 0203 Armenia
E-mail: aishkhanyan@gmail.com

We examine the power-series solutions and the series solutions in terms of the Hermite functions for the biconfluent Heun equation. Infinitely many cases for which a solution of the biconfluent equation is presented as an irreducible linear combination of four generalized hypergeometric functions, that in general do not reduce to polynomials, are identified.

**MSC numbers**: 33E30 Other functions coming from differential, difference and integral equations, 34B30 Special equations (Mathieu, Hill, Bessel, etc.), 30Bxx Series expansions

**Keywords:** Bi-confluent Heun equation, series expansion, recurrence relation

## 1. Introduction

Over the recent years, the five functions of the Heun class have been finding extensive applications in various branches of physics and mathematics (see, e.g., [1,2] and references therein). These functions arise when one solves certain types of linear second-order ordinary differential equations that present natural generalizations of the equations satisfied by the hypergeometric functions. At present time, the most complete information on the Heun functions is presented in [3-5].

In the present paper, we examine the finite-sum reductions of the power and Hermite-function series expansions of the biconfluent Heun functions. The result we report is that there are an infinitely many cases for which the biconfluent Heun equation permits a solution as irreducible linear combination of four generalized hypergeometric functions.

## 2. Power-series solution of the biconfluent Heun equation

Consider the biconfluent Heun equation [3-5] which we write in the following form:

$$z\frac{d^2\Phi}{dz^2} + \left(p_0 + p_1 s z - 2s^2 z^2\right)\frac{d\Phi}{dz} + \left(q_0 s + q_1 s^2 z\right)\Phi = 0, \quad (1)$$

where $p_{0,1,2}$, $q_{0,1}$ and $s$ are arbitrary complex constants. Though the number of the independent parameters involved in this equation can readily be reduced to four by scaling and shifting the argument $z$, however, this form is more convenient for our treatment. We note that the accessory parameter $q_0$ is the only one which does not characterize the behavior of the solution in the vicinity of a singular point (see [3,4]).

The coefficients of the Frobenius power-series expansion [6] of the particular solution in the vicinity of the regular singular point $z = 0$ corresponding to the zero characteristic exponent, that is the coefficients $c_n$ of a solution

$$\Phi(z) = \sum_{n=0}^{\infty} c_n (sz)^n \qquad (2)$$

(where the factor $s$ is put for convenience), obey the three-term recurrence relation [3-5]

$$R_n c_n + Q_{n-1} c_{n-1} + P_{n-2} c_{n-2} = 0, \qquad (3)$$

with

$$R_n = n(n-1+p_0), \quad Q_n = p_1 n + q_0, \quad P_n = -2n + q_1. \qquad (4)$$

In explicit form, the recurrence relation reads

$$n(n-1+p_0)c_n + \left(p_1(n-1) + q_0\right)c_{n-1} + \left(-2(n-2) + q_1\right)c_{n-2} = 0. \qquad (5)$$

The solution (2) becomes a polynomial of degree $N$ if $P_N = 0$, i.e.

$$q_1 = 2N, \quad N = 0, 1, 2..., \qquad (6)$$

and $-q_0$ is an eigenvalue of the $(N+1) \times (N+1)$ minor $M^{(N)}$ of the infinite matrix

$$M = \begin{bmatrix} 0 & p_0 & 0 & 0 & \dots \\ 2N & p_1 & 2(1+p_0) & 0 & \dots \\ 0 & 2(N-1) & 2p_1 & 3(2+p_0) & \dots \\ \dots & \dots & \dots & \dots & \dots \end{bmatrix}. \qquad (7)$$

It is readily understood that the eigenvalue equation indeed presents a $(N+1)$ degree polynomial equation for the accessory parameter $q_0$. For $N = 1, 2, 3$ these equations are explicitly written as

$N = 0$, $q_1 = 0$: $\qquad\qquad\qquad q_0 = 0,$ \hfill (8)

$N = 1$, $q_1 = 2$: $\qquad\qquad\qquad q_0^2 + p_1 q_0 - 2p_0 = 0,$ \hfill (9)

$N = 2$, $q_1 = 4$: $\qquad q_0^3 + 3p_1 q_0^2 + (2p_1^2 - 8p_0 - 4)q_0 - 8p_0 p_1 = 0,$ \hfill (10)

$N = 3$, $q_1 = 6$:
$$q_0^4 + 6p_1 q_0^3 + (11p_1^2 - 20p_0 - 20)q_0^2 + 6(-10p_0 p_1 - 6p_1 + p_1^3)q_0 + (-36p_0 p_1^2 + 36p_0^2 + 72p_0) = 0. \qquad (11)$$

Note that for a given eigenvalue $-q_0$ the expansion coefficients $c_0$, ..., $c_N$ are the components of the eigenvector of the corresponding minor $M^{(N)}$.

## 3. Expansion of the biconfluent Heun function in terms of the Hermite functions

Under rather general assumptions, a solution of equation (1) can be presented in the form of a series in terms of the Hermite functions of a scaled and shifted argument [7,8]:

$$\Phi(z) = \sum_{n=0}^{\infty} d_n H_{p_0+q_1/2+n}(sz - p_1/2), \qquad (12)$$

where the coefficients $d_n$ are calculated using the recurrence relation (3) with coefficients

$$R_n = 2n(n + p_0 + q_1/2), \quad Q_n = (n + p_0)p_1 - q_0, \quad P_n = n + p_0, \qquad (13)$$

so that the recurrence relation explicitly reads

$$2n(n + p_0 + q_1/2)d_n + (p_1(n-1) + p_0 p_1 - q_0)d_{n-1} + (n - 2 + p_0)d_{n-2} = 0. \qquad (14)$$

The series (12) is reduced to a finite sum, for which the maximum number of nonzero coefficients is equal to $N+1$, if

$$p_0 = -N, \quad N = 0, 1, 2..., \qquad (15)$$

and $q_0$ is the eigenvalue of the $(N+1)\times(N+1)$ minor of the matrix $M$ with elements

$$M = \begin{bmatrix} -Np_1 & 2(1-N)+q_1 & 0 & 0 & \cdots \\ -N & (1-N)p_1 & 4(2-N)+2q_1 & 0 & \cdots \\ 0 & 1-N & (2-N)p_1 & 6(3-N)+3q_1 & \cdots \\ \cdots & \cdots & \cdots & \cdots & \cdots \end{bmatrix}. \qquad (16)$$

Again, the eigenvalue equation presents a $(N+1)$-degree polynomial equation for $q_0$. The first four eigenvalue equations are explicitly written as

$N = 0$, $p_0 = 0$: $\qquad\qquad\qquad q_0 = 0,$ $\qquad\qquad(17)$

$N = 1$, $p_0 = -1$: $\qquad\qquad q_0^2 + p_1 q_0 + q_1 = 0,$ $\qquad\qquad(18)$

$N = 2$, $p_0 = -2$: $\quad q_0^3 + 3p_1 q_0^2 + (2p_1^2 + 4q_1 - 4)q_0 + 4p_1 q_1 = 0,$ $\qquad(19)$

$N = 3$, $p_0 = -3$: 
$$q_0^4 + 6p_1 q_0^3 + (10q_1 + 11p_1^2 - 20)q_0^2 + \\ (30p_1 q_1 + 6p_1^3 - 36p_1)q_0 + 9q_1^2 + (18p_1^2 - 36)q_1 = 0. \qquad (20)$$

As in the case of the above power-series expansion, the coefficients $d_0$, ..., $d_N$ of expansion (12) are components of the eigenvector of the corresponding minor of matrix (16).

## 4. Special power series solutions of the biconfluent Heun equation

An observation concerning the power series expansion (2)-(5) for the vicinity of the singularity $z = 0$ is that if the involved parameters are such that the series in terms of the Hermite functions terminates into a finite sum (that is if $p_0$ is a non-positive integer, see equation (15)), the condition $R_n \neq 0$ for the power series (2)-(5) is violated at $n = 1 - p_0$ (see equation (4)). Hence, in this case equations (4) cannot be applied to construct a power-series expansion with the zero characteristic exponent (the other Frobenius solution with the greater exponent $\mu = 1 - p_0$ [3,4] is of course applicable).

A further observation is that in the case of non-positive integer $p_0 = 0, -1, -2, ...$ one can derive useful results by trying the power-series expansions of the biconfluent Heun function $\Phi(z)$ in the vicinity of a *regular* point $z = z_0 \neq 0$. In constructing such expansions, it turns out that it is rather helpful to start from the representation of the function $\Phi(z)$ as a linear combination of the Hermite functions and further apply the power-series expansions of the involved Hermite functions, for instance, the following ones [5]:

$$H_\nu(z) = \sum_{k=0}^{\infty} \frac{\Gamma(k/2 - \nu/2)}{2\Gamma(-\nu)} \cdot \frac{(-2z)^k}{k!}, \tag{21}$$

$$H_\nu(z + \xi) = \sum_{k=0}^{\infty} (-\nu)_k \cdot H_{\nu-k}(\xi) \cdot \frac{(-2z)^k}{k!}. \tag{22}$$

We note that as a result of application of equation (22) to the expansion point $z_0 = 0$ one arrives at an expansion of function $\Phi(z)$ in powers of $z$ defined by the coefficients

$$c_k = \sum_{j=0}^{N} d_j \cdot (-p_0 - q_1/2 - j)_k \cdot H_{p_0 + q_1/2 + j - k}(-p_1/2) \frac{(-2)^k}{k!}. \tag{23}$$

An observation now is that for these coefficients no two-term recurrence relation is known.

Notably, the expansion of function $\Phi(z)$ in the vicinity of the point $z_0 = p_1/(2s)$, in other words, the expansion in powers of the *scaled and shifted* variable $\xi = sz - p_1/2$:

$$\Phi(z) = \sum_{n=0}^{\infty} \tilde{c}_n \xi^n, \tag{24}$$

leads to two-term recurrence relations for certain infinite subsets of expansion coefficients. Note that in this case the expansion coefficients are explicitly written as

$$\tilde{c}_k = \sum_{j=0}^{N} d_j \frac{\Gamma(k/2 - p_0/2 - q_1/4 - j/2)}{2\Gamma(-p_0 - q_1/2 - j)} \frac{(-2)^k}{k!}. \tag{25}$$

Using this formula, we will show that if the Hermite-function expansion (12) is terminated, the solution of the biconfluent Heun equation can be written in terms of the generalized-hypergeometric functions.

To proceed to the derivation of the latter result, let

$$\beta = p_0 + q_1/2, \tag{26}$$

$N'$ is the integer part of $N/2$ and $N''$ is the integer part of $(N-1)/2$.

Now, as the first step, one can show that the coefficients $\tilde{c}_k$ can be decomposed as

$$\tilde{c}_k = c'_k + c''_k, \tag{27}$$

where $c'_k$ and $c''_k$ are separately written in a form similar to (25) using either even or odd values of the summation index $j$:

$$c'_{2k} = \frac{\Gamma(-\beta/2-N')}{2\Gamma(-\beta-2N')} \cdot \frac{(-\beta/2-N')_k \cdot (-2)^{2k}}{(2k)!} \sum_{j=0}^{N'} \frac{(-\beta/2-N'+k)_j}{(-\beta-2N')_{2j}} \cdot d_{2N'-2j}, \tag{28}$$

$$c'_{2k+1} = \frac{\Gamma(1/2-\beta/2-N')}{2\Gamma(-\beta-2N')} \frac{(1/2-\beta/2-N')_k(-2)^{2k+1}}{(2k+1)!} \sum_{j=0}^{N'} \frac{(1/2-\beta/2-N'+k)_j}{(-\beta-2N')_{2j}} \cdot d_{2N'-2j}, \tag{29}$$

$$c''_{2k} = \frac{\Gamma(-1/2-\beta/2-N'')}{2\Gamma(-1-\beta-2N'')} \cdot \frac{(-1/2-\beta/2-N'')_k \cdot (-2)^{2k}}{(2k)!}$$
$$\sum_{j=0}^{N''} \frac{(-1/2-\beta/2-N''+k)_j}{(-1-\beta-2N'')_{2j}} \cdot d_{2N''-2j+1}, \tag{30}$$

$$c''_{2k+1} = \frac{\Gamma(-\beta/2-N'')}{2\Gamma(-1-\beta-2N')} \cdot \frac{(-\beta/2-N'')_k \cdot (-2)^{2k+1}}{(2k+1)!} \sum_{j=0}^{N''} \frac{(-\beta/2-N''+k)_j}{(-1-\beta-2N'')_{2j}} \cdot d_{2N''-2j+1}, \tag{31}$$

It is next shown that these coefficients obey the following two-term recurrence relations:

$$\frac{c'_{2k+2}}{c'_{2k}} = \frac{(-\beta/2-N'+k)}{(k+1)(k+1/2)} \cdot \frac{\varphi(k+1)}{\varphi(k)}, \tag{32}$$

$$\frac{c'_{2k+3}}{c'_{2k+1}} = \frac{(1/2-\beta/2-N'+k)}{(k+1)(k+3/2)} \cdot \frac{\varphi(k+3/2)}{\varphi(k+1/2)}, \tag{33}$$

$$\frac{c''_{2k+2}}{c''_{2k}} = \frac{(-1/2-\beta/2-N''+k)}{(k+1)(k+1/2)} \cdot \frac{\psi(k+1/2)}{\psi(k-1/2)}, \tag{34}$$

$$\frac{c''_{2k+3}}{c''_{2k+1}} = \frac{(-\beta/2-N''+k)}{(k+1)(k+3/2)} \cdot \frac{\psi(k+3/2)}{\psi(k+1/2)}, \tag{35}$$

where

$$\varphi(z) = \sum_{j=0}^{N'} \frac{(-\beta/2-N'+z)_j}{(-\beta-2N')_{2j}} \cdot d_{2N'-2j} \tag{36}$$

is a polynomial of degree $N'$ and

$$\psi(z) = \sum_{j=0}^{N''} \frac{(-\beta/2 - N'' + z)_j}{(-1 - \beta - 2N'')_{2j}} \cdot d_{2N''-2j+1} \qquad (37)$$

is a polynomial of degree $N''$.

Note that these recurrence relations are applicable with the proviso that

$$\tilde{c}_0 = \sum_{j=0}^{N} \frac{d_j \cdot \sqrt{\pi} \cdot 2^{\beta+j}}{\Gamma(1/2 - \beta/2 - j/2)}, \qquad (38)$$

$$\tilde{c}_1 = -\sum_{j=0}^{N} \frac{d_j \cdot \sqrt{\pi} \cdot 2^{1+\beta+j}}{\Gamma(-\beta/2 - j/2)}, \qquad (39)$$

$$c'_0 = \frac{\sqrt{\pi} \, 2^{\beta}}{\Gamma(1/2 - \beta/2)} \sum_{j=0}^{N'} (-4)^j \cdot d_{2j} \cdot (1/2 + \beta/2)_j, \qquad (40)$$

$$c'_1 = -\frac{\sqrt{\pi} \, 2^{1+\beta}}{\Gamma(-\beta/2)} \sum_{j=0}^{N'} (-4)^j \cdot d_{2j} \cdot (1 + \beta/2)_j, \qquad (41)$$

$$c''_0 = \frac{\sqrt{\pi} \, 2^{1+\beta}}{\Gamma(-\beta/2)} \sum_{j=0}^{N''} (-4)^j \cdot d_{2j+1} \cdot (1 + \beta/2)_j, \qquad (42)$$

$$c''_1 = -\frac{\sqrt{\pi} \, 2^{2+\beta}}{\Gamma(-1/2 - \beta/2)} \sum_{j=0}^{N''} (-4)^j \cdot d_{2j+1} \cdot (3/2 + \beta/2)_j. \qquad (43)$$

## 5. Auxiliary relations for generalized hypergeometric functions

To proceed further, we need some particular relations concerning the generalized hypergeometric function

$$F(z) = {}_pF_q(\alpha_1, \alpha_2, ..., \alpha_p; \gamma_1, ..., \gamma_q; z). \qquad (44)$$

Using the differentiation formula for hypergeometric functions [5,9,10]

$$\frac{d^n}{dz^n}\left(z^{\alpha_1+n-1}F(\omega z)\right) = (\alpha_1)_n \cdot z^{\alpha_1-1} \cdot {}_pF_q(\alpha_1 + n, \alpha_2, ..., \alpha_p; \gamma_1, ..., \gamma_q; \omega z), \qquad (45)$$

one can show that

$${}_pF_q(\alpha_1 + n, \alpha_2, ..., \alpha_p; \gamma_1, ..., \gamma_q; z) = \sum_{k=0}^{n} \frac{C_n^k}{(\alpha_1)_k} \cdot z^k \left(\frac{d}{dz}\right)^k F(z), \qquad (46)$$

where $C_n^k = \binom{n}{k}$ are binomial coefficients. Furthermore, it can be shown by mathematical induction that

$$z^k \left(\frac{d}{dz}\right)^k = \left(z\frac{d}{dz}\right)\left(z\frac{d}{dz}-1\right)\cdots\left(z\frac{d}{dz}-k+1\right). \tag{47}$$

With this, an arbitrary linear combination of generalized hypergeometric functions $_pF_q(\alpha_1+n,\alpha_2,...,\alpha_p;\gamma_1,...,\gamma_q;z)$, $n=0,1,2,...,N$, can be presented in the following form:

$$\sum_{n=0}^{N} b_n\, _pF_q(\alpha_1+n,\alpha_2,...,\alpha_p;\gamma_1,...,\gamma_q;z)$$

$$= \sum_{n=0}^{N}\sum_{k=0}^{n} b_n \frac{C_n^k}{(\alpha_1)_k}\left(z\frac{d}{dz}\right)\left(z\frac{d}{dz}-1\right)\cdots\left(z\frac{d}{dz}-k+1\right)F(z)$$

$$= \frac{b_N}{(\alpha_1)_N}\left(z\frac{d}{dz}+\lambda_1\right)\cdots\left(z\frac{d}{dz}+\lambda_N\right)F(z), \tag{48}$$

where $\lambda_1,\ ...,\ \lambda_N$ are the roots multiplied by $(-1)$ of the $N$ th degree polynomial

$$g(\xi) = \sum_{n=0}^{N}\sum_{k=0}^{n} b_n \frac{C_n^k}{(\alpha_1)_k} \xi(\xi-1)\cdot...\cdot(\xi-k+1). \tag{49}$$

The first four terms of the double sum involved in this equation are explicitly written as

$$g(\xi) = b_0 + b_1\left(1+\frac{\xi}{\alpha_1}\right) + b_2\left(1+\frac{2\xi}{\alpha_1}+\frac{\xi(\xi-1)}{\alpha_1(\alpha_1+1)}\right)$$

$$+ b_3\left(1+\frac{3\xi}{\alpha_1}+\frac{3\xi(\xi-1)}{\alpha_1(\alpha_1+1)}+\frac{\xi(\xi-1)(\xi-2)}{\alpha_1(\alpha_1+1)(\alpha_1+2)}\right)+.... \tag{50}$$

Note that the coefficient of the highest-degree term $\xi^N$ of $g(\xi)$ is equal to $b_N/(\alpha_1)_N$ and the free term is equal to

$$\frac{b_N}{(\alpha_1)_N}\cdot\lambda_1\cdot...\cdot\lambda_N = \sum_{k=0}^{N} b_k. \tag{51}$$

Using another formula for differentiation of hypergeometric functions [10]:

$$\frac{d}{dz}\left(z^\lambda\cdot F(\omega z)\right) = \lambda\cdot z^{\lambda-1}\cdot\,_{p+1}F_{q+1}(\alpha_1,...,\alpha_p,\lambda+1;\gamma_1,...,\gamma_q,\lambda;\omega z), \tag{52}$$

another useful relation is derived:

$$\lambda\,_{p+1}F_{q+1}(\alpha_1,...,\alpha_p,\lambda+1;\gamma_1,...,\gamma_q,\lambda;z) = z^{-\lambda}\cdot z\frac{d}{dz}\left(z^\lambda F(z)\right) = \left(z\frac{d}{dz}+\lambda\right)F(z). \tag{53}$$

By applying this formula for several times, we then have

$$\lambda_1\cdot...\cdot\lambda_N\cdot\,_{p+N}F_{q+N}(\alpha_1,...,\alpha_p,\lambda_1+1,...,\lambda_N+1;\gamma_1,...,\gamma_q,\lambda_1,...,\lambda_N;z) =$$

$$\left(z\frac{d}{dz}+\lambda_1\right)...\left(z\frac{d}{dz}+\lambda_N\right)F(z). \tag{54}$$

As a result, one arrives at the relation

$$\sum_{n=0}^{N} b_n \, _pF_q(\alpha_1+n,\alpha_2,...,\alpha_p;\gamma_1,...,\gamma_q;z) =$$

$$\frac{b_N}{(\alpha_1)_N}\lambda_1\cdot...\cdot\lambda_N \cdot {}_{p+N}F_{q+N}(\alpha_1,...,\alpha_p,\lambda_1+1,...,\lambda_N+1;\gamma_1,...,\gamma_q,\lambda_1,...,\lambda_N;z) =$$

$$(b_0+...+b_N)\cdot {}_{p+N}F_{q+N}(\alpha_1,...,\alpha_p,\lambda_1+1,...,\lambda_N+1;\gamma_1,...,\gamma_q,\lambda_1,...,\lambda_N;z). \tag{55}$$

In particular, for $N=1$ we have:

$$b_0 \, _pF_q(\alpha_1,\alpha_2,...,\alpha_p;\gamma_1,...,\gamma_q;z) + b_1 \, _pF_q(\alpha_1+1,\alpha_2,...,\alpha_p;\gamma_1,...,\gamma_q;z) =$$

$$(b_0+b_1)\cdot {}_{p+1}F_{q+1}(\alpha_1,...,\alpha_p,\lambda+1;\gamma_1,...,\gamma_q,\lambda;z), \tag{56}$$

where

$$\lambda = \alpha_1(1+b_0/b_1). \tag{57}$$

Similarly, for $N=2$ the result reads

$$\sum_{n=0}^{2} b_n \, _pF_q(\alpha_1+n,\alpha_2,...,\alpha_p;\gamma_1,...,\gamma_q;z) =$$

$$(b_0+b_1+b_2)\cdot {}_{p+1}F_{q+1}(\alpha_1,...,\alpha_p,\lambda_1+1,\lambda_2+1;\gamma_1,...,\gamma_q,\lambda_1,\lambda_2;z), \tag{58}$$

where $\lambda_1$ and $\lambda_2$ are the roots of the quadratic equation

$$b_2\frac{\lambda(\lambda+1)}{\alpha_1(\alpha_1+1)} - (b_1+2b_2)\frac{\lambda}{\alpha_1} + (b_0+b_1+b_2) = 0. \tag{59}$$

## 6. Generalized-hypergeometric solutions of the biconfluent Heun equation

We now proceed to the derivation of generalized-hypergeometric solutions of the biconfluent Heun equation (1). Let

$$\beta = p_0 + q_1/2, \quad \xi = sz - p_1/2, \tag{60}$$

$N'$ is the integer part of $N/2$ and $N''$ is the integer part of $(N-1)/2$. The Hermite functions in terms of which the expansion (12) is developed have the following representation in terms of the Kummer confluent hypergeometric functions [5]:

$$u_n(z) = H_{\beta+n}(\xi) = \frac{2^{\beta+n}\sqrt{\pi}}{\Gamma(1/2-\beta/2-n/2)}{}_1F_1(-\beta/2-n/2;1/2;\xi^2) -$$

$$\frac{2^{\beta+n+1}\sqrt{\pi}}{\Gamma(-\beta/2-n/2)}\xi\,_1F_1(1/2-\beta/2-n/2;3/2;\xi^2). \tag{61}$$

With this, the solution (12) is rewritten as

$$\Phi(z) = \sum_{n=0}^{N'} d_{2n} \cdot \frac{2^{\beta+2n}\sqrt{\pi}}{\Gamma(1/2-\beta/2-n)}\,_1F_1(-\beta/2-n;1/2;\xi^2)$$

$$-\sum_{n=0}^{N'} d_{2n} \cdot \frac{2^{\beta+1+2n}\sqrt{\pi}}{\Gamma(-\beta/2-n)}\xi\,_1F_1(1/2-\beta/2-n;3/2;\xi^2)$$

$$+\sum_{n=0}^{N''} d_{2n+1} \cdot \frac{2^{\beta+1+2n}\sqrt{\pi}}{\Gamma(-\beta/2-n)}\,_1F_1(-\beta/2-1/2-n;1/2;\xi^2)$$

$$-\sum_{n=0}^{N''} d_{2n+1} \cdot \frac{2^{\beta+2+2n}\sqrt{\pi}}{\Gamma(-\beta/2-1/2-n)}\xi\,_1F_1(-\beta/2-n;3/2;\xi^2), \tag{62}$$

or otherwise

$$\Phi(z) = \frac{2^{\beta+2N'}\sqrt{\pi}}{\Gamma(1/2-\beta/2-N')}\sum_{n=0}^{N'} \frac{d_{2N'-2n}}{4^n(1/2-\beta/2-N')_n}\,_1F_1(-\beta/2-N'+n;1/2;\xi^2)$$

$$-\frac{2^{\beta+1+2N'}\sqrt{\pi}}{\Gamma(-\beta/2-N')}\xi\sum_{n=0}^{N'} \frac{d_{2N'-2n}}{4^n(-\beta/2-N')_n}\,_1F_1(1/2-\beta/2-N'+n;3/2;\xi^2)$$

$$+\frac{2^{\beta+1+2N''}\sqrt{\pi}}{\Gamma(-\beta/2-N'')}\sum_{n=0}^{N''} \frac{d_{2N''-2n+1}}{4^n(-\beta/2-N'')_n}\,_1F_1(-\beta/2-1/2-N''+n;1/2;\xi^2)$$

$$-\frac{2^{\beta+2+2N''}\sqrt{\pi}}{\Gamma(-\beta/2-1/2-N'')}\xi\sum_{n=0}^{N''} \frac{d_{2N''-2n+1}}{4^n(-\beta/2-1/2-N'')_n}\,_1F_1(-\beta/2-N''+n;3/2;\xi^2). \tag{63}$$

Using the above formula (55) we then obtain

$$\Phi(z) = C' \cdot {}_{1+N'}F_{1+N'}(-\beta/2-N',\lambda'_1+1,...,\lambda'_{N'}+1;1/2,\lambda'_1,...,\lambda'_{N'};\xi^2)$$

$$-D'\,\xi \cdot {}_{1+N'}F_{1+N'}(1/2-\beta/2-N',\mu'_1+1,...,\mu'_{N'}+1;3/2,\mu'_1,...,\mu'_{N'};\xi^2)$$

$$+C'' \cdot {}_{1+N''}F_{1+N''}(-\beta/2-1/2-N'',\lambda''_1+1,...,\lambda''_{N'}+1;1/2,\lambda''_1,...,\lambda''_{N'};\xi^2)$$

$$-D''\,\xi \cdot {}_{1+N''}F_{1+N''}(-\beta/2-N'',\mu''_1+1,...,\mu''_{N'}+1;3/2,\mu''_1,...,\mu''_{N'};\xi^2), \tag{64}$$

where

$$C' = \frac{2^{\beta+2N'}\sqrt{\pi}}{\Gamma(1/2-\beta/2-N')}\sum_{k=0}^{N'} \frac{d_{2N'-2n}}{4^n(1/2-\beta/2-N')_n}, \tag{65}$$

$$D' = \frac{2^{\beta+1+2N'}\sqrt{\pi}}{\Gamma(-\beta/2-N')}\sum_{k=0}^{N'}\frac{d_{2N'-2n}}{4^n(-\beta/2-N')_n}, \quad (66)$$

$$C'' = \frac{2^{\beta+1+2N''}\sqrt{\pi}}{\Gamma(-\beta/2-N'')}\sum_{k=0}^{N''}\frac{d_{2N''-2n+1}}{4^n(-\beta/2-N'')_n}, \quad (67)$$

$$D'' = \frac{2^{\beta+2+2N''}\sqrt{\pi}}{\Gamma(-\beta/2-1/2-N'')}\sum_{k=0}^{N''}\frac{d_{2N''-2n+1}}{4^n(-\beta/2-1/2-N'')_n} \quad (68)$$

with $\lambda'_1, \ldots, \lambda'_{N'}$ being the roots of the polynomial

$$g'(\xi) = \sum_{n=0}^{N'}\sum_{k=0}^{n}\frac{d_{2N'-2n}}{4^n(1/2-\beta/2-N')_n}\cdot\frac{C_n^k}{(-\beta/2-N')_k}\cdot\xi(\xi-1)\cdot\ldots\cdot(\xi-k+1), \quad (69)$$

$\mu'_1, \ldots, \mu'_{N'}$ being the roots of the polynomial

$$h'(\xi) = \sum_{n=0}^{N'}\sum_{k=0}^{n}\frac{d_{2N'-2n}}{4^n(-\beta/2-N')_n}\cdot\frac{C_n^k}{(1/2-\beta/2-N')_k}\cdot\xi(\xi-1)\cdot\ldots\cdot(\xi-k+1), \quad (70)$$

$\lambda''_1, \ldots, \lambda''_{N''}$ being the roots of polynomial

$$g''(\xi) = \sum_{n=0}^{N''}\sum_{k=0}^{n}\frac{d_{2N''-2n+1}}{4^n(-\beta/2-N'')_n}\cdot\frac{C_n^k}{(-\beta/2-1/2-N'')_k}\cdot\xi(\xi-1)\cdot\ldots\cdot(\xi-k+1), \quad (71)$$

and $\mu''_1, \ldots, \mu''_{N''}$ being the roots of polynomial

$$h''(\xi) = \sum_{n=0}^{N''}\sum_{k=0}^{n}\frac{d_{2N''-2n+1}}{4^n(-\beta/2-1/2-N'')_n}\cdot\frac{C_n^k}{(-\beta/2-N'')_k}\cdot\xi(\xi-1)\cdot\ldots\cdot(\xi-k+1). \quad (72)$$

For the first three values of $N$ we have the following explicit solutions.

$N = 0$:

$$\Phi(z) = \frac{2^{\beta}\sqrt{\pi}}{\Gamma(1/2-\beta/2)}{}_1F_1(-\beta/2;1/2;\xi^2) - \frac{2^{\beta+1}\sqrt{\pi}}{\Gamma(-\beta/2)}\xi\,{}_1F_1(1/2-\beta/2;3/2;\xi^2). \quad (73)$$

$N = 1$:

$$\Phi(z) = d_0\frac{2^{\beta}\sqrt{\pi}}{\Gamma(1/2-\beta/2)}{}_1F_1(-\beta/2;1/2;\xi^2) + d_1\frac{2^{\beta+1}\sqrt{\pi}}{\Gamma(-\beta/2)}{}_1F_1(-\beta/2-1/2;1/2;\xi^2)$$

$$-d_0\frac{2^{\beta+1}\sqrt{\pi}}{\Gamma(-\beta/2)}\xi\,{}_1F_1(1/2-\beta/2;3/2;\xi^2) - d_1\frac{2^{\beta+2}\sqrt{\pi}}{\Gamma(-\beta/2-1/2)}\xi\,{}_1F_1(-\beta/2;3/2;\xi^2), \quad (74)$$

where
$$d_0 = 1, \quad d_1 = \frac{p_1+q_0}{2(1+\beta)}, \quad (75)$$

$N = 2$:

$$\Phi(z) = \frac{2^{\beta}\sqrt{\pi}}{\Gamma(-\beta/2+1/2)} C \ _2F_2(-\beta/2-1,(-\beta/2-1)C+1;1/2,(-\beta/2-1)C;\xi^2)$$

$$+d_1\frac{2^{\beta+1}\sqrt{\pi}}{\Gamma(-\beta/2)} {}_1F_1(-\beta/2-1/2;1/2;\xi^2)$$

$$-\frac{2^{\beta+1}\sqrt{\pi}}{\Gamma(-\beta/2)} D\,\xi\, {}_2F_2(-\beta/2-1/2,(-\beta/2-1/2)D+1;3/2,(-\beta/2-1/2)D;\xi^2)$$

$$-d_1\frac{2^{\beta+2}\sqrt{\pi}}{\Gamma(-\beta/2-1/2)} \xi\, {}_1F_1(-\beta/2;3/2;\xi^2), \qquad (76)$$

where
$$d_0 = 1, \quad d_1 = \frac{2p_1+q_0}{2(1+\beta)}, \quad d_2 = \frac{2+(p_1+q_0)d_1}{4(2+\beta)}, \qquad (77)$$

and
$$C = 1+4d_2(-\beta/2-1/2), \qquad D = 1+4d_2(-\beta/2-1). \qquad (78)$$

$N = 3$:

$$\Phi(z) = \frac{2^{\beta}\sqrt{\pi}}{\Gamma(-\beta/2+1/2)} \cdot d_0\, C'\, {}_2F_2(-\beta/2-1,(-\beta/2-1)C'+1;1/2,(-\beta/2-1)C';\xi^2)$$

$$+\frac{2^{\beta+1}\sqrt{\pi}}{\Gamma(-\beta/2)} \cdot d_1\, C''\, {}_2F_2(-\beta/2-3/2,(-\beta/2-3/2)C''+1;1/2,(-\beta/2-3/2)C'';\xi^2)$$

$$-\frac{2^{\beta+1}\sqrt{\pi}}{\Gamma(-\beta/2)} \cdot d_0\, D'\,\xi\, {}_2F_2(-\beta/2-1/2,(-\beta/2-1/2)D'+1;3/2,(-\beta/2-1/2)D';\xi^2)$$

$$-\frac{2^{\beta+2}\sqrt{\pi}}{\Gamma(-\beta/2-1/2)} \cdot d_1\, D''\,\xi\, {}_2F_2(-\beta/2-1,(-\beta/2-1)D''+1;3/2,(-\beta/2-1)D'';\xi^2), \quad (79)$$

where
$$d_0 = 1, \quad d_1 = \frac{3p_1+q_0}{2(1+\beta)}, \quad d_2 = \frac{3+(2p_1+q_0)d_1}{4(2+\beta)}, \quad d_3 = \frac{2d_1+(p_1+q_0)d_2}{6(3+\beta)}, \qquad (80)$$

and
$$C' = 1+4(-\beta/2-1/2)d_2/d_0, \qquad D' = 1+4(-\beta/2-1)d_2/d_0, \qquad (81)$$

$$C'' = 1+4(-\beta/2-1)d_3/d_1, \qquad D'' = 1+4(-\beta/2-3/2)d_3/d_1. \qquad (82)$$

The solution (64) with the parameters (65)-(72) presents the main result of the present paper.

## 7. Discussion

Series solutions of the Heun equations have been discussed by many authors. In particular, expansions in terms of the hypergeometric functions of various types have been constructed (see, e.g., [11-22]). For these series, the finite-sum reductions play an important role since the hypergeometric functions are well studied.

Generalized-hypergeometric solutions of the Heun equations have been first reported by Lettessier who derived a few solutions for the general and single-confluent Heun equations [23-25], see also [26,27]. Recently, these results have been generalized to show that there exist an infinite number of such solutions for the two mentioned Heun equations [28-30].

In the present paper, we have shown that generalized-hypergeometric solutions are also possible for the biconfluent Heun equation. Though these solutions are in a sense similar to those obtained earlier for the general and confluent Heun equations, however, there exist essential peculiarities. The main difference is that in the biconfluent case the solutions are presented as irreducible linear combinations of *four* generalized-hypergeometric functions, while in the general and single-confluent cases the solution is written through a single generalized-hypergeometric function. Another difference is that in the latter cases coefficients of the three-term recurrence relations governing the power-series expansions are reduced to two-term ones, while for the biconfluent Heun equation this property is not the case.

**Acknowledgments**

This research was supported by the Russian-Armenian (Slavonic) University at the expense of the Ministry of Education and Science of the Russian Federation, the Armenian Science Committee (SC Grants No. 18RF-139 and 18T-1C276), and the Armenian National Science and Education Fund (ANSEF Grant No. PS-4986).